\documentclass[11pt]{article}
\usepackage{amsmath}
\usepackage{amssymb}
\usepackage{amscd}
\usepackage[all]{xy}
\usepackage[T1]{fontenc}
\usepackage{epsfig}
\def\qsl{\mathcal O(SL_q(2))}
\def\fre{k[z,z^{-1}] * \mathcal O(SL_q(2))}
\def\N{\mathbb N}
\def\Z{\mathbb Z}

\newtheorem{theo}{Theorem}[section]
\newtheorem{prop}[theo]{Proposition}
\newtheorem{coro}[theo]{Corollary}
\newtheorem{lemm}[theo]{Lemma}
\newtheorem{defi}[theo]{Definition}

\title{\textbf{\textsc{corepresentation theory of universal cosovereign hopf algebras}}} 
\author{\textsc{Julien Bichon}}
\date{{\small \textsl{Laboratoire de Math\'ematiques Appliqu\'ees,
Universit\'e de Pau et des Pays de l'Adour, \\
IPRA, Avenue de l'universit\'e,
64000 Pau, France.}
E-mail: Julien.Bichon@univ-pau.fr}}

\makeatletter
\renewcommand{\@makefnmark}{}
\makeatother

\begin{document}

\maketitle

\begin{abstract}
We determine the corepresentation theory of the universal
cosovereign Hopf algebras, for generic matrices over an algebraically
closed field of characteristic zero. Our results generalize
Banica's previous results in the compact case.
As an application, we easily get the representation
theory of the quantum automorphism group
of a matrix algebra endowed with a non-necessarily tracial measure. 
\end{abstract}

\section{Introduction}

The universal cosovereign Hopf algebras are  natural free analogues
of the general linear groups in quantum group theory.
In this paper we determine their corepresentation theory in the
generic case.

Let $k$ be commutative field and let $F \in GL(n,k)$.
The algebra $H(F)$ \cite{[Bi2]} is defined to be the universal algebra
with generators
 $(u_{ij})_{1 \leq i,j \leq n}$,
 $(v_{ij})_{1 \leq i,j \leq n}$ and relations:
$$ {u} {^t \! v} = {^t \! v} u = I_n \quad ; \quad {vF} {^t \! u} F^{-1} = 
{F} {^t \! u} F^{-1}v = I_n,$$
where $u= (u_{ij})$, $v = (v_{ij})$ and $I_n$ is
the identity $n \times n$ matrix. It turns out \cite{[Bi2]} that
$H(F)$ is a Hopf algebra
with comultiplication defined by 
$\Delta(u_{ij}) = \sum_k u_{ik} \otimes u_{kj}$
and $\Delta(v_{ij}) = \sum_k v_{ik} \otimes v_{kj}$, with counit
defined by $\varepsilon (u_{ij}) = \varepsilon (v_{ij}) = \delta_{ij}$ and
with antipode defined by $S(u) = {^t \! v}$ and $S(v) = F {^t \! u} F^{-1}$.
Furthermore $H(F)$ is a cosovereign Hopf algebra \cite{[Bi2]}: there
exists an algebra morphism $\Phi : H(F) \longrightarrow k$
such that $S^2 = \Phi * {\rm id} * \Phi^{-1}$. The Hopf algebras 
$H(F)$ have the following universal property (\cite{[Bi2]}, Theorem 3.2).

\smallskip

\textsl{Let $A$ be Hopf algebra and let $V$ be finite dimensional
$A$-comodule isomorphic with its bidual comodule $V^{**}$. 
Then there exists a matrix $F \in GL(n,k)$ ($n = \dim V$) such that
$V$ is an $H(F)$-comodule and such that there exists a Hopf algebra
morphism $\pi : H(F) \longrightarrow A$ such that 
$(1_V \otimes \pi) \circ \beta_V = \alpha_V$,
where $\alpha_V : V \longrightarrow V \otimes A$ and
$\beta_V : V \longrightarrow V \otimes H(F)$ denote the coactions 
of $A$ and $H(F)$ on $V$ respectively. In particular 
every finite type cosovereign Hopf algebra is a homomorphic quotient
of a Hopf algebra $H(F)$.
}

\smallskip

In view of this universal property it is natural to say that the 
Hopf algebras $H(F)$ are the universal cosovereign Hopf algebras,
or the free cosovereign Hopf algebras, and to see these Hopf algebras
as natural analogues of the general linear groups in quantum
group theory. It also seems to be important to describe their
corepresentation theory.

\smallskip

If $k = \mathbb C$ and if $F$ is a positive matrix, the Hopf algebras
$H(F)$ are nothing but the CQG algebras associated to the
universal compact quantum groups introduced by Van Daele and Wang
\cite{[VDW]}.
In this case the corepresentation theory has been worked out by 
Banica \cite{[Ba]}: the simple corepresentations correspond to the 
elements of the free product $\mathbb N * \mathbb N$, and the fusion rules
are described by an ingenious formula involving a new product
$\odot$ on the free algebra on two generators.
We generalize Banica's results to the case of the case of an arbitrary 
generic matrix, over any algebraically closed field of characteristic zero.

\medskip

In order to state our main result, we need to introduce some notation
and terminology.  

\noindent
$\bullet$
Let $F \in GL(n,k)$. We say that $F$ is normalized if ${\rm tr}(F) =
{\rm tr} (F^{-1})$. We say that $F$ is normalizable if there exists 
$\lambda \in k^*$ such that 
${\rm tr}(\lambda F) = {\rm tr} ((\lambda F)^{-1})$. Over an algebraically 
closed field, any matrix is normalizable unless 
${\rm tr}(F) = 0 \not = {\rm tr} (F^{-1})$  or
${\rm tr}(F)\not = 0 = {\rm tr} (F^{-1})$.
We will only essentially consider normalized matrices $F$ or equivalently
normalizable matrices, since $H(\lambda F) = H(F)$.

\noindent
$\bullet$ Let $q$ in $k^*$. As usual, we say that $q$ is generic
if $q$ is not a root of unity of order $N \geq 3$.
We say that a matrix $F \in GL(n,k)$ is generic if $F$ is normalized
and if the solutions of $q^2 -{\rm tr}(F)q +1 = 0$ are generic. 

\noindent 
$\bullet$ 
Let $q$ in $k^*$. We put $F_q = \left(\begin{array}{cc} q^{-1} & 0 \\
                          0 & q \\
       \end{array} \right) \in GL(2,k)$.
The Hopf algebra $H(F_q)$ is denoted by $H(q)$.

\noindent
$\bullet$ Let $F \in GL(n,k)$. The natural $n$-dimensional 
$H(F)$-comodules associated to the multiplicative matrices $u = (u_{ij})$
and $v = (v_{ij})$ are denoted by $U$ and $V$, with $V = U^*$.

\noindent
$\bullet$ We will consider the coproduct monoid $\mathbb N * \mathbb N$.
Equivalently $\mathbb N * \mathbb N$ is the free monoid on two generators,
which we denote, as in \cite{[Ba]}, by $\alpha$ and $\beta$.
There is a unique antimultiplicative morphism
$^- : \mathbb N * \mathbb N \longrightarrow \mathbb N * \mathbb N$
such that $\bar{e} = e$, $\bar{\alpha} = \beta$ and $\bar{\beta} = \alpha$
($e$ denotes the unit element of $\mathbb N * \mathbb N$).

\medskip

We can now state our main result. Here $k$ denotes an algebraically
closed field.

\begin{theo}
Let $F \in GL(n,k)$ ($n \geq 2$) be a normalized matrix.

\noindent
a) Let  $q \in k^*$ be such that $q^2 -{\rm tr}(F)q +1 = 0$. The 
comodule categories over $H(F)$ and $H(q)$ are monoidally equivalent.

We assume now that $k$ is a characteristic zero field.

\noindent
b) The Hopf algebra $H(F)$ is cosemisimple if and only if $F$ is 
a generic matrix.

\noindent
c) Assume that $F$ is generic. 
To any element $x \in \mathbb N * \mathbb N$ corresponds a simple
$H(F)$-comodule $U_x$, with $U_e = k$, $U_\alpha = U$ and $U_\beta = V$.
Any simple $H(F)$-comodule is isomorphic with one of the $U_x$, and
$U_x \cong U_y$ if and only if $x=y$. For $x,y \in \mathbb N * \mathbb N$, 
we have $U_x^* \cong U_{\bar{x}}$ and
$$U_x \otimes U_y  \cong \bigoplus_{\{a,b,g \in 
\mathbb N * \mathbb N|x=ag,y={\bar g}b\}} U_{ab} \ .$$ 
\end{theo}

It is clear from the statement that the proof of Theorem 1.1 is
divided into two parts. The first part reduces the corepresentation theory
of $H(F)$ to the one of $H(q)$. Then we realize $H(q)$ as a Hopf
subalgebra of the free product $\fre$, and we conclude using the 
classification of simple comodules of a free product of cosemisimple 
Hopf algebras \cite{[Wa1]}, and Banica's product
$\odot$ on the free algebra on two generators.	 

\medskip

Another interesting class of universal Hopf algebras 
was constructed by Wang \cite{[Wa2]} in the compact quantum
group framework: these are the quantum automorphism groups of 
finite-dimensional (measured) $C^*$-algebras. 
The corepresentation theory,
similar to the one of $SO(3)$, was described by Banica \cite{[Ba2]},
for $C^*$-algebras endowed with (good) tracial measures.
A special case 
of a general construction of \cite{[Bi1]} yields algebraic analogues
of Wang's quantum automorphism groups. Using the previous results
concerning $H(F)$,
it is not  difficult to describe the
representation theory
of the quantum automorphism group of a matrix algebra endowed with
a non-necessarily tracial measure, reducing the computations to the
case of the quantum group $SO_q(3)$.

\medskip

This paper is organized as follows. In Section 2 we use the Hopf-Galois
systems techniques of \cite{[Bi4]} to show that for a normalized matrix
$F$, there exists $q \in  k^*$ such that the comodule categories over 
$H(F)$ and $H(q)$ are monoidally equivalent. 
This section also includes results for non-normalizable matrices.
In Section 3 we construct an
injective algebra morphism of $H(q)$ into the free product algebra
$\fre$. In Section 4 we show that $H(q)$ is cosemisimple if and only if
$q$ is generic, and Section 5 contains the classification of the
simple $H(q)$-comodules and their fusion rules in the generic case.
Section 6 is devoted to some applications of Theorem 1.1 to 
structural properties of the Hopf algebras $H(F)$. Finally in Section 7
we use our previous results to describe the representation category
of the quantum automorphism group of a matrix algebra endowed with
a non-necessarily tracial measure. 

\section{Reduction to the two-dimensional case}

This section is essentially devoted to prove part a) of Theorem 1.1.
In fact we consider a more general situation and get results
for non-normalizable matrices.
We will use Hopf-Galois systems techniques \cite{[Bi4]}.
We will not repeat here the definition a Hopf-Galois system,
for which we refer to \cite{[Bi4]}.

\medskip 

Let $E \in GL(m,k)$ and let $F \in GL(n,k)$. Recall \cite{[Bi4]} that
the algebra $H(E,F)$ is the universal algebra with generators
$u_{ij}$, $ v_{ij}$, $1\leq i \leq m, 1\leq j \leq n$,
and satisfying the relations
$$u {^t \! v} = I_m = v F {^t \! u} E^{-1} \quad ; \quad
{^t \! v}u = I_n = F {^t \! u} E^{-1} v.$$
When $E = F$, we have $H(F,F) = H(F)$. In fact the
Hopf algebra structure of $H(F)$ is just a particular case
of the fact that if $H(E,F)$ is a non-zero algebra, then 
$(H(E),H(F),\break H(E,F),H(F,E))$ is a Hopf-Galois system (see Proposition
4.3 in \cite{[Bi4]}). Combining Proposition 4.3 and Corollary 1.4
in \cite{[Bi4]}, we have the following result.

\begin{prop}
Assume that $H(E,F)$ is a non-zero algebra. Then the comodule categories
over $H(E)$ and $H(F)$ are monoidally equivalent. $\square$
\end{prop}

So we have to study the algebras $H(E,F)$. 
It is not difficult to see that if $H(E,F) \not = \{0 \}$, then 
${\rm tr}(E) = {\rm tr}(F)$ and ${\rm tr}(E^{-1}) = {\rm tr}(F^{-1})$. 
The converse assertion will essentially follow from the next result, where
some technical conditions are required.

\begin{prop}
Let $E \in GL(m,k)$ and let $F \in GL(n,k)$ ($m,n \geq 2$).
Assume that $E$ is a diagonal matrix, that $F$ is a lower-triangular
matrix, that 
${\rm tr}(E) = {\rm tr}(F)$ and ${\rm tr}(E^{-1}) = {\rm tr}(F^{-1})$. 
Then the elements ($u_{ij}$), $1 \leq i \leq  m$, $1\leq j \leq n$, 
generate
a free subalgebra on $mn$ generators. In particular $H(E,F)$ is a non zero-algebra.
\end{prop} 

As in \cite{[Bi3]}, we will use the diamond lemma \cite{[Be]}. 
Let us write down
explicitly a presentation of $H(E,F)$: $H(E,F)$ is the universal
algebra with generators $u_{ij}$, $ v_{ij}$, $1\leq i \leq m, 1\leq j \leq n$,
and relations:
\begin{align}
& u_{in}v_{jn} = \delta_{ij} - \sum_{k=1}^{n-1}u_{ik}v_{jk} \ , \quad
1 \leq i,j \leq m \\
&  v_{i1}u_{j1} = F_{11}^{-1}(E_{ij} - \sum_{k=2, l=1}^{n}F_{kl}
v_{ik}u_{jl}) \ , \quad
1 \leq i,j \leq m \\ 
&  v_{1i}u_{1j} = \delta_{ij} - \sum_{k=2}^{m}v_{ki}u_{kj} \ , \quad
1 \leq i,j \leq n \\
& u_{mi}v_{mj} = E_{mm}(F_{ij}^{-1} - \sum_{k=1}^{m-1}
E_{kk}^{-1}u_{ki}v_{kj}) \ , \quad
1 \leq i,j \leq n.  
\end{align}
We have a nice presentation to use the diamond lemma \cite{[Be]}, 
of which we freely use the techniques and definitions.
We only need the simplified exposition of \cite{[KS]}.
We order the set of monomials in the following way.
We order the set $\{1, \ldots, m \} \times \{1, \ldots, n \}$
lexicographically. Then we order the set $\{u_{ij} \}$ with the order induced 
by the preceding order, and we order the set $\{v_{ij}\}$ with the
inverse order. We order the set $X= \{u_{ij}, v_{kl}\}$ in such a way
that $v_{11} < u_{11}$. Finally two monomials are ordered according
to their length, and two monomials of equal length are ordered 
lexicographically according to the order on the set $X$.
It is clear that the order just defined is compatible
with the above presentation. 

\begin{lemm}
There are exactly two inclusion ambiguities: $(v_{11}u_{11},v_{11}u_{11})$
and $(u_{mn}v_{mn}, \break u_{mn}v_{mn})$. There are exactly the following
overlap ambiguities.
$$(u_{in}v_{1n}, v_{1n}u_{1j}) , \quad 
(v_{i1}u_{m1}, u_{m1}v_{mj}) ,  \quad 1\leq i \leq m, \quad
1 \leq j \leq n.$$ 
$$ (v_{1i}u_{1n}, u_{1n}v_{jn}) , \quad
(u_{mi}v_{m1}, v_{m1}u_{j1}) , \quad 
1\leq i \leq n, \quad 1 \leq j \leq m.$$ 
All these ambiguities are resolvable.
\end{lemm}

\noindent
\textbf{Proof}. It is easy to see that the ambiguities above are the only
ones. 
Let us check that the first inclusion ambiguity is resolvable.
 As usual the
symbol ``$\rightarrow$'' means that we perform a reduction.
We have:
\begin{align*}
& F_{11}^{-1}(E_{11} - \sum_{k=2,l=1}^n F_{kl}v_{1k}u_{1l})
\rightarrow F_{11}^{-1}(E_{11} - \sum_{k=2,l=1}^n F_{kl}
(\delta_{kl}-\sum_{r=2}^m v_{rk} u_{rl})) \\
& = F_{11}^{-1}(E_{11} - \sum_{k=2}^n F_{kk} +
\sum_{k=2,l=1}^n \sum_{r=2}^m F_{kl}v_{rk}u_{rl})
\end{align*}
On the other hand we have:
\begin{align*}
& 1 - \sum_{k=2}^m v_{k1}u_{k1} \rightarrow
1- \sum_{k=2}^mF_{11}^{-1}(E_{kk} - \sum_{l=2,r=1}^n
F_{lr}v_{kl}u_{kr}) \\
& = F_{11}^{-1}(E_{11} - \sum_{k=2}^n F_{kk} +
\sum_{k=2,l=1}^n \sum_{r=2}^m F_{kl}v_{rk}u_{rl})
\end{align*}
We have used the identity ${\rm tr}(E) = {\rm tr}(F)$.
Hence this inclusion ambiguity is resolvable.
Also it is not difficult to check,
using ${\rm tr}(E^{-1}) = {\rm tr}(F^{-1})$,
that the other inclusion ambiguity is resolvable.
This is left to the reader.
Let us check that the first two families of overlap
ambiguities are resolvable. The resolvability of
the other two families will be left to the reader. 
First consider $(u_{in}v_{1n}, v_{1n}u_{1j}),  1\leq i \leq m,
1 \leq j \leq n.$ We have:
$$
(\delta_{i1} - \sum_{k=1}^{n-1}u_{ik} v_{1k})u_{1j}
\rightarrow \delta_{i1} u_{1j} - \sum_{k=1}^{n-1}
u_{ik}(\delta_{kj} -\sum_{l=2}^m v_{lk}u_{lj})
 = \delta_{i1} u_{1j} - (1- \delta_{jn})u_{ij} +
\sum_{k=1}^{n-1} \sum_{l=2}^m u_{ik} v_{lk} u_{lj}.$$
On the other hand we have:
\begin{align*}
& u_{in}(\delta_{nj} - \sum_{k=2}^{m}v_{kn} u_{kj}) \rightarrow
\delta_{nj}u_{in} -\sum_{k=2}^m(\delta_{ki}- \sum_{l=1}^{n-1}
u_{il}v_{kl})u_{kj} = \\ 
& \delta_{nj}u_{in} - (1-\delta_{1i})u_{ij} + 
\sum_{k=2}^m \sum_{l=1}^{n-1} u_{il}v_{kl}u_{kj}
= \delta_{i1} u_{1j} - (1- \delta_{jn})u_{ij} +
\sum_{k=1}^{n-1} \sum_{l=2}^m u_{ik} v_{lk} u_{lj}.
\end{align*}
Hence these ambiguities are resolvable. Let us now study the ambiguities
$(v_{i1}u_{m1}, u_{m1}v_{mj})$, $1\leq i \leq m, 1 \leq j \leq n.$
We have:
\begin{align*}
& F_{11}^{-1}(E_{im} - \sum_{k=2,l=1}^n F_{kl}v_{ik}u_{ml})v_{mj}\\
& \rightarrow
F_{11}^{-1}(E_{im}v_{mj} - \sum_{k=2,l=1}^n(F_{kl}v_{ik}(
E_{mm}(F_{lj}^{-1}- \sum_{r=1}^{m-1}E_{rr}^{-1}u_{rl}v_{rj}))) \\
& = F_{11}^{-1}E_{mm}(\delta_{im}v_{mj} - (1-\delta_{j1}) v_{ij}
+ \sum_{k=2,l=1}^n \sum_{r=1}^{m-1} F_{kl}E_{rr}^{-1}
v_{ik}u_{rl}v_{rj}).
\end{align*}
On the other hand we have:
\begin{align*}
& v_{i1}(E_{mm}(F_{1j}^{-1}-\sum_{k=1}^{m-1}E_{kk}^{-1}u_{k1}v_{kj}) \\
& \rightarrow E_{mm}(\delta_{1j}F_{11}^{-1}v_{i1} -
\sum_{k=1}^{m-1}E_{kk}^{-1}(F_{11}^{-1}(\delta_{ik} E_{ii}
- \sum_{l=2,r=1}^n F_{lr}v_{il}u_{kr}))v_{kj}) \\
& = E_{mm} F_{11}^{-1}(\delta_{1j} v_{i1} - (1 - \delta_{im})v_{ij}
+ \sum_{k=1}^{m-1} \sum_{l=2,r=1}^{n-1}E_{kk}^{-1} F_{lr}
v_{il}u_{kr}v_{kj}) \\
& = F_{11}^{-1}E_{mm}(\delta_{im}v_{mj} - (1-\delta_{j1}) v_{ij}
+ \sum_{k=2,l=1}^n \sum_{r=1}^{m-1} F_{kl}E_{rr}^{-1}
v_{ik}u_{rl}v_{rj}).
\end{align*}
Hence these ambiguities are resolvable. $\square$

\bigskip

\noindent
\textbf{Proof of Proposition 2.2}. Since our order is compatible
with the presentation, and since all the ambiguities are resolvable, we can 
use the diamond lemma \cite{[Be]}: the reduced monomials form
a basis of $H(E,F)$, and in particular the monomials in 
elements of the set $\{u_{ij}, 1 \leq i \leq m, 1 \leq j \leq n \}$
are linearly independent, and hence the elements of this set
generate a free subalgebra on $mn$ generators. In particular 
$H(E,F)$ is a non-zero algebra. $\square$

\bigskip

We can now easily prove the following slightly more general result.

\begin{prop} 
Let $E \in GL(m,k)$ and let $F \in GL(n,k)$ ($m,n \geq 2$).
Assume that
${\rm tr}(E) = {\rm tr}(F)$ and ${\rm tr}(E^{-1}) = {\rm tr}(F^{-1})$. 
Then $H(E,F)$ is a non-zero algebra.
\end{prop}

\noindent
\textbf{Proof}. Since we want to prove that $H(E,F)$ is a non-zero
algebra, we can assume that $k$ is algebraically closed.
For matrices $P \in GL(m,k)$ and let $Q \in GL(n,k)$, the algebras
$H(E,F)$ and $H(PEP^{-1},QFQ^{-1})$ are isomorphic (\cite{[Bi4]},
Proposition 4.2), thus we can assume that the matrices $E$ and $F$ are
lower-triangular. Consider $G \in GL(m,k)$ a diagonal
matrix such that   
${\rm tr}(G) = {\rm tr}(E) = {\rm tr}(F)$ and 
${\rm tr}(G^{-1}) = {\rm tr}(E^{-1}) = {\rm tr}(F^{-1})$.
By the proof of Proposition 4.3 in \cite{[Bi4]}, there exists an algebra 
morphism
$\delta : H(E,F) \longrightarrow H(E,G) \otimes H(G,F)$ such
that $\delta(u_{ij}) =  \sum_{k=1}^m
u_{ik} \otimes u_{kj}$. Also there exists an algebra morphism
$\phi : H(E,G) \longrightarrow H(G,E)^{\rm op}$ such that 
$\phi(u) = {^t \! v}$. Thus we have an algebra morphism
 $\delta ': H(E,F) \longrightarrow H(G,E)^{\rm op} \otimes H(G,F)$ such
that $\delta '(u_{ij}) =  \sum_{k=1}^m
v_{ki} \otimes u_{kj}$. By the proof of Proposition 2.2, the elements
$(v_{ij})$, $(u_{ij})$ are linearly independent elements of
$H(G,E)$ and $H(G,F)$ respectively. Hence it is clear that
$H(E,F)$ is a non-zero algebra. $\square$

\bigskip

Combining Propositions 2.1 and 2.4, we can state the main result of
the section, which contains part a) of Theorem 1.1. Recall that
for $q \in k^*$, we put $H(q) = H(F_q)$ where 
$F_q = \left(\begin{array}{cc} q^{-1} & 0 \\
                          0 & q \\
       \end{array} \right) \in GL(2,k)$.
Note also that if $k$ is algebraically closed,
any $F \in GL(2,k)$ is normalizable.

\begin{coro} Let $F \in GL(n,k)$ ($n \geq 2$)
and assume that $k$ is algebraically closed.

\noindent
a) Assume that $F$ is normalizable. Then there exists $q \in k^*$ such
that we have an equivalence of monoidal categories:
$${\rm Comod} (H(F)) \cong^\otimes {\rm Comod}(H(q)).$$
If $F$ is normalized, we take $q$ as a solution of the equation
$q^2 -{\rm tr}(F)q +1 = 0$.

\noindent
b) Assume that $F$ is not normalizable. Let $E \in GL(3,k)$ be any matrix
such that ${\rm tr}(E) = 0$ and ${\rm tr}(E^{-1}) \not = 0$. Then we have
an equivalence of monoidal categories:
$${\rm Comod} (H(F)) \cong^\otimes {\rm Comod}(H(E)).$$
\end{coro}

\noindent
\textbf{Proof}. a) Let $\lambda \in k^*$ be such that
${\rm tr}(\lambda F) = {\rm tr} ((\lambda F)^{-1})$, and let
$q \in k^*$ be a solution of
$q^2 -{\rm tr}(\lambda F)q +1 = 0$. This equation is equivalent to
${\rm tr}(F_q^{-1}) =
{\rm tr}(F_q) = q + q^{-1} = 
{\rm tr}(\lambda F) = {\rm tr}((\lambda F)^{-1}))$. 
By Proposition 2.4, $H(F_q,F)$ is a non-zero algebra, and 
we conclude using Proposition 2.1.

\noindent
b) Since $F$ is not normalizable and since the base field is algebraically 
closed, we have ${\rm tr}(F) = 0 \not = {\rm tr} (F^{-1})$  or
${\rm tr}(F)\not = 0 = {\rm tr} (F^{-1})$. Since the Hopf
algebras $H(F)$ and $H(^t \!F^{-1})$ are isomorphic
(\cite{[Bi2]}, Proposition 3.3), we can assume
that ${\rm tr}(F) = 0 \not = {\rm tr} (F^{-1})$. Since $k$
is algebraically closed, there always exists 
$E \in GL(3,k)$ satisfying ${\rm tr}(E) = 0$ and ${\rm tr}(E^{-1}) \not = 0$,
and we conclude as in part a). $\square$  

\bigskip

Recall that the fundamental $n$-dimensional comodule of $H(F)$ 
associated to the multiplicative matrix $(u_{ij})$ is denoted
by $U$. The following result reflects the ``freeness'' of $H(F)$.

\begin{coro}
 Let $F \in GL(n,k)$. The comodules $U^{\otimes k}$, $k \in \N$, are
simple non-equivalent $H(F)$-comodules.
\end{coro}

\noindent
\textbf{Proof}. We can assume that $k$ is algebraically closed.
If $n=1$ then $H(F)$ is just the algebra of Laurent polynomials
$k[z,z^{-1}]$, so the result is immediate. Assume now that
$n \geq 2$. First assume that $F$ is a diagonal matrix.
By Proposition 2.2 the monomials in the elements
$u_{ij}$ form a linearly independent subset of $H(F)$, and hence
the comodules $U^{\otimes k}$, $k \in \N$, are
simple non-equivalent $H(F)$-comodules. Now assume that
$F$ is a lower-triangular matrix. Take $E\in GL(n,k)$ a diagonal matrix such 
that ${\rm tr}(E) = {\rm tr}(F)$ and ${\rm tr}(E^{-1})
={\rm tr}(F^{-1})$. The monoidal category equivalence of Proposition 2.1
transforms the $H(F)$-comodule $U$ into the $H(E)$-comodule $U$
(see \cite{[U],[Sc1],[Bi4]} for the construction). Hence we conclude 
by the diagonal case. This finishes the proof since 
the Hopf algebras $H(PFP^{-1})$ and $H(F)$ are isomorphic 
for $P \in GL(m,k)$ (\cite{[Bi2]}). $\square$ 

\begin{coro}
 Let $F \in GL(n,k)$ be a non-normalizable matrix.
Then the Hopf algebra $H(F)$ is not cosemisimple.
\end{coro}

\noindent
\textbf{Proof}. By the preceding corollary $U$ is a simple $H(F)$-comodule. 
We can assume that $k$ is algebraically closed.
If $H(F)$ was cosemisimple, and since
${^t \! F}^{-1}$ is an intertwiner between $U$ and $U^{**}$
(see the Proof of Theorem 3.2 in \cite{[Bi2]}), then we would
have by (\cite{[KS]}, Proposition 15, chapter 11,
 or the original reference \cite{[La]}) ${\rm tr}(F) \not = 0$
and ${\rm tr}(F^{-1}) \not = 0$, which would contradict our assumption.
$\square$

\section{The algebra $H(q)$}

This section is devoted to the construction of an algebra
embedding of $H(q) = H(F_q)$ into $\fre$.
This embedding will be used later to study the corepresentation
theory of $H(q)$.

\medskip

Let $q \in k^*$. The algebra $H(q)$ has 8 generators.
We put $\alpha = u_{11}$, $\beta = u_{12}$, 
$\gamma = u_{21}$, $\delta = u_{22}$, 
$\alpha^* = v_{11}$, $\beta^*= v_{12}$, 
$\gamma^* = v_{21}$, $\delta^* = v_{22}$. 
Let us rewrite the presentation of $H(q)$:
it is the universal algebra with generators 
$\alpha$, $\beta$, $\gamma$, $\delta$, $\alpha^*$, $\beta^*$, 
$\gamma^*$, $\delta^*$ and satisfying the relations:
$$\left\lbrace
\begin{array}{l}
\beta  \beta^* = 1 - \alpha \alpha^* \\
\beta \delta^* = - \alpha \gamma^* \\
\delta \beta^* = - \gamma \alpha^* \\
\delta \delta^* = 1 - \gamma \gamma^*
\end{array}
\right.
\left\lbrace
\begin{array}{l}
\alpha^*  \alpha = 1 - q^2\beta^* \beta \\
\alpha^* \gamma = - q ^2 \beta^* \delta \\
\gamma^* \alpha = - q^2 \delta^* \beta \\
\gamma^* \gamma = q^2(1 - \delta^* \delta)
\end{array}
\right.
\left\lbrace
\begin{array}{l}
\alpha^*  \alpha = 1 - \gamma^* \gamma \\
\alpha^* \beta = - \gamma^* \delta \\
\beta^* \alpha = - \delta^* \gamma \\
\beta^* \beta = 1 - \delta^* \delta
\end{array}
\right.
\left\lbrace
\begin{array}{l}
\gamma  \gamma^* = q^2(1 - \alpha \alpha^*) \\
\gamma \delta^* = - q^2 \alpha \beta^* \\
\delta \gamma^* = - q^2 \beta \alpha^* \\
\delta \delta^* = 1 - q^2 \beta \beta^*
\end{array}
\right.$$
Note that the fourth relation of the first family and that the
first relation of the second family are redundant.
We have left these redundant relations in order to use
the results of Section 2, where some redundant relations were
also present.

\medskip

We define now an algebra extension of $H(q)$, which will
be denoted by $H^+(q)$. This algebra will be shown to be isomorphic with
$\fre$.

\begin{defi}
The algebra $H^+(q)$ is the universal algebra with generators 
$\alpha$, $\beta$, $\gamma$, $\delta$, $\alpha^*$, $\beta^*$, 
$\gamma^*$, $\delta^*$, $t$, $t^{-1}$,
and satisfying the relations of $H(q)$ and:
$$t t^{-1} = 1 = t^{-1}t \ ; \ t^{-1} \alpha = \delta^{*} t \ ; \
t^{-1} \beta = -q^{-1}\gamma^{*} t \ ; \ 
t^{-1} \gamma  = -q\beta^{*} t \ ; \
t^{-1} \delta = \alpha^{*} t.$$
\end{defi}

There is an obvious algebra morphism $H(q) \longrightarrow H^+(q)$.

\begin{lemm}
The natural algebra morphism $H(q) \longrightarrow H^+(q)$ is
injective.
\end{lemm}

\noindent
\textbf{Proof}. We will use again the diamond lemma, since we have not been
able to find a more direct way to prove our lemma.
First we order the set 
$\{\alpha, \beta, \gamma, \delta, \alpha^*, \beta^*, 
\gamma^*, \delta^*, t, t^{-1}\}$ in the following way:
$$\delta^* < \gamma^* < \beta^* < \alpha^* < 
\alpha < \beta < \gamma < \delta< t^{-1} < t.$$
Two monomials of different length are ordered according to their length
and two monomials of equal length are ordered lexicographically
according to the above order. In order to resolve some ambiguities,
let us rewrite the presentation of $H^+(q)$:
$H^+(q)$ is the universal algebra with generators 
$\alpha$, $\beta$, $\gamma$, $\delta$, $\alpha^*$, $\beta^*$, 
$\gamma^*$, $\delta^*$, $t$, $t^{-1}$,
and satisfying the relations of $H(q)$ and
\begin{align*}
& tt^{-1} = t^{-1}t \ ; \ t^{-1}t = 1 \ ; \
t^{-1} \alpha = \delta^{*} t \ ; \ t \delta^* = \alpha t^{-1} \ ; \
t^{-1} \beta = -q^{-1}\gamma^* t \ ; \ t \gamma^* = -q \beta t^{-1} \ ; \\ 
& t^{-1} \gamma  = -q\beta^* t \ ; \ t \beta^* = -q^{-1} \gamma t^{-1} \ ; \
t^{-1} \delta = \alpha^{*} t \ ; \ t \alpha^* = \delta t^{-1}.
\end{align*}
It is clear the order just defined is compatible with this
presentation. There are the ambiguities of Lemma 2.3, which  were
shown to be resolvable there, there are no other inclusion ambiguities
and the following overlap ambiguities:
\begin{align*}
& (tt^{-1},t^{-1}\alpha) \quad ; \quad  (tt^{-1},t^{-1}\beta) \quad ; \quad
(tt^{-1},t^{-1}\gamma) \quad ; \quad (tt^{-1},t^{-1}\delta) \quad ; \\
& (t^{-1}t,t \delta^*) \quad ; \quad (t^{-1}t,t \gamma^*) \quad ; \quad
(t^{-1}t,t \beta^*) \quad ; \quad (t^{-1}t,t \alpha^*) \quad ; \\
& (t^{-1}\beta, \beta \beta^*) \quad ; \quad (t^{-1}\beta, \beta \delta^*) 
\quad ; \quad
(t\gamma^*, \gamma^* \alpha) \quad ; \quad (t\gamma^*, \gamma^* \gamma) 
\quad ; \\
& (t^{-1}\gamma, \gamma \gamma^*) \quad ; \quad 
(t^{-1}\gamma, \gamma \delta^*) \quad ; \quad
(t\beta^*, \beta^* \alpha) \quad ; \quad (t\beta^*, \beta^* \beta) \quad ; \\
& (t^{-1}\delta, \delta \beta^*) \quad ; \quad 
(t^{-1}\delta, \delta \delta^*) \quad ; \quad
(t^{-1}\delta, \delta \gamma^*) \quad ; \quad 
(t^{-1}\delta, \delta \delta^*) \quad ; \\
& (t \alpha^*, \alpha^* \alpha) \quad ; \quad 
(t \alpha^*, \alpha^* \gamma) \quad ; \quad
(t \alpha^*, \alpha^* \alpha) \quad ; \quad (t\alpha^*, \alpha^* \beta). 
\end{align*}  
These ambiguities are easily seen to be resolvable: this is left to the reader.
Hence by the diamond lemma the reduced monomials form a basis of
$H^+(q)$. It is clear that the reduced monomials of $H(q)$
(for the reductions of Section 2) are still reduced monomials
in $H^+(q)$, and hence the images under $H(q) \rightarrow H^+(q)$
of the elements of a basis of $H(q)$ are still linealy independant
elements, which proves that our algebra map is injective. $\square$

\bigskip

Recall that $\qsl$ is the universal algebra with generators
$a,b,c,d$ and relations
$$ba = qab \ ; \ ca = qac \ ; \ db = qbd \ ; \ 
dc = qcd \ ; \ cb = bc = q(ad-1) \ ; da = qbc +1.$$
The algebra just defined is $\mathcal O(SL_{q^{-1}}(2))$ in 
\cite{[KS]}. Our convention does not change the resulting Hopf
algebra, up to isomorphism.
Now consider the free product $\fre$, that is the coproduct of 
$k[z,z^{-1}]$ and of $\qsl$ in the category of unital algebras.
We have the following result:

\begin{lemm}
There exists a unique algebra isomorphism
$\tilde{\pi}: H^+(q) \rightarrow \fre$
such that
\begin{align*}
& \tilde{\pi}(\alpha) = za, \ 
\tilde{\pi}(\beta) = zb, \ \tilde{\pi}(\gamma) = zc, \ 
\tilde{\pi}(\delta) = zd, \
\tilde{\pi}(\alpha^*) = dz^{-1}, \
\tilde{\pi}(\beta^*) = -q^{-1}cz^{-1}, \\ 
& \tilde{\pi}(\gamma^*)= -qbz^{-1}, \ \tilde{\pi}(\delta^*) =az^{-1},
\ \tilde{\pi}(t) = z, \ \tilde{\pi}(t^{-1}) = z^{-1}.
\end{align*}  
\end{lemm}

\noindent
\textbf{Proof}. It is a direct verification to check the
existence of the algebra morphism $\tilde{\pi}$.
Let us construct an inverse isomorphism. First there is an algebra
morphism $\rho_1 : k[z,z^{-1}] \longrightarrow H^+(q)$ defined by 
$\rho_1(z) = t$. It is also a direct verification to check
the existence of an algebra morphism $\rho_2 :
\qsl \longrightarrow H^+(q)$ such that
$$\rho_2(a) = t^{-1} \alpha = \delta^* t, \
\rho_2(b) = t^{-1} \beta = -q^{-1} \gamma^* t, \
\rho_2(c) = t^{-1} \gamma = -q \beta^*t, \
\rho_2(d) = t^{-1} \delta = \alpha^* t.$$
Using the universal property of the free product,
we have a unique algebra morphism $\rho :
\fre \longrightarrow H^+(q)$ extending $\rho_1$ and $\rho_2$.
It is straightforward to check that $\tilde \pi$ and $\rho$ are
mutually inverse isomorphisms. $\square$

\bigskip

We arrive at the main result of the section.

\begin{prop}
There exists an injective algebra morphism
$\pi : H(q) \longrightarrow \fre$ such that
\begin{align*}
& \pi(\alpha) = za \ , \ 
\pi(\beta) = zb \ , \ \pi(\gamma) = zc \ , \ 
\pi(\delta) = zd \ , \\ 
& \pi(\alpha^*) = dz^{-1}, \
\pi(\beta^*) = -q^{-1}cz^{-1},  
\pi(\gamma^*)= -qbz^{-1}, \ \pi(\delta^*) =az^{-1}.
\end{align*}
\end{prop}
 
\noindent
\textbf{Proof}. The algebra morphism announced is just the composition
of the injective algebra morphisms of Lemmas 3.2 and 3.3,
so is itself injective $\square$

\section{Cosemisimplicity of $H(q)$}

In this section, where $k$ is assumed to be an algebraically closed
field of characteristic zero, we show that $H(q)$ is cosemisimple
if and only if $q$ is generic.

\medskip

First let us recall that if $A$ and $B$ are Hopf algebras,
their free product may be endowed with a natural Hopf algebra
structure, induced by the Hopf algebras structures of $A$ and $B$.
For example $\fre$ is a Hopf algebra, and by a straightforward 
verification, we have the following result.

\begin{prop}
The injective algebra morphism
$\pi : H(q) \longrightarrow \fre$
is a Hopf algebra morphism. $\square$
\end{prop}

Wang \cite{[Wa1]} has studied free products of Hopf algebras at
the compact quantum group level.
His results may be adapted to arbitrary cosemisimple Hopf algebras 
without difficulties.
Let us recall the main results. In the following
$A$ and $B$ denote cosemisimple Hopf algebras.

\noindent
$\bullet$ The Hopf algebra $A * B$ is still cosemisimple.
This may be shown as follows. Consider the Haar functionals
(see e.g. \cite{[KS]}) $h_A$ and $h_B$ on $A$ and $B$ respectively,
and form their free product $h_A*h_B$ as in \cite{[Av]}, Proposition 1.1.
Then $h_A*h_B$ is a Haar functional on $A*B$ (see \cite{[Wa1]}, Theorem 3.8)
and thus $A*B$ is a cosemisimple Hopf algebra

\smallskip

\noindent
$\bullet$ An $A*B$-comodule is said to be a simple alternated
$A*B$-comodule if it has the form $V_1 \otimes \ldots \otimes V_n$,
where each $V_i$ is a simple non-trivial $A$-comodule or $B$-comodule, and
if $V_i$ is an $A$-comodule, then $V_{i+1}$ is a $B$-comodule, and
conversely. A simple alternated $A*B$-comodule is a simple $A*B$-comodule,
and every non-trivial simple $A*B$-comodule is isomorphic
with a simple alternated $A*B$-comodule (see \cite{[Wa1]}, Theorem 3.10).  

\smallskip

\noindent
$\bullet$ Let $V$ and $W$ be simple alternated $A*B$-comodules.
Assume that $V$ ends by an $A$-comodule and that $W$ begins by
a $B$-comodule. Then $V \otimes W$ is decomposed into a direct sum
of simple alternated comodules according to the decomposition
of tensor products of $A$-comodules. The same thing holds for $B$.  

\medskip

We will use these results to prove the following fact.

\begin{prop}
Let $q \in k^*$. Then $H(q)$ is cosemisimple if and only
if $q$ is generic.
\end{prop}

\noindent
\textbf{Proof}. We will use the following well-known fact.
Let $A \subset B$ be a Hopf algebra inclusion. Then an $A$-comodule
is semisimple if and only if it is semisimple as a $B$-comodule.
In particular if $B$ is cosemisimple, so is $A$.
First assume that $q$ is generic. Then it is well-known that
$\qsl$ is cosemisimple (see e.g. \cite{[KS]}), and since
$k[z,z^{-1}]$ is also cosemisimple, we have that $\fre$ is cosemisimple,
and so is $H(q)$ by Proposition 4.1.

Let us now assume that $q$ is a root of unity of order $N \geq 3$.
We will construct a non semisimple $H(q)$-comodule.
Put $N_0 = N/2$ if $N$ is even and $N_0 = N$ if $N$ is odd. Let 
$V_1$ be the fundamental two-dimensional $\qsl$-comodule. One 
can deduce from the results of \cite{[KP]} that 
$V_1^{\otimes N_0}$ is not a semisimple $\qsl$-comodule. For $i \in \Z$ 
we denote by $Z^i$ the one-dimensional comodule associated to the group-like
element $z^i$ of $k[z,z^{-1}]$. Using $\pi$, we view $H(q)$ as a Hopf 
subalgebra of $\fre$ and by the construction of $\pi$, $Z \otimes V_1$
and $V_1 \otimes Z^{-1}$ are $H(q)$-comodules. Then 
$V_1^{\otimes 2} = V_1 \otimes Z^{-1} \otimes Z \otimes V_1$ is an 
$H(q)$-comodule. Assume that $N_0$ is even: $N_0 = 2k$. Then
$V_1^{\otimes N_0} = V_1^{\otimes 2k}$ is an $H(q)$-comodule. Since
$V_1^{\otimes N_0}$ is not a semisimple $\qsl$-comodule, it is not
a semisimple $\fre$-comodule, and so is not a semisimple 
$H(q)$-comodule. Assume now that $N_0$ is odd: $N_0 = 2k+1$. We have
seen that $V_1^{\otimes 2k}$ is an $H(q)$-comodule, and hence 
$Z \otimes V_1^{\otimes N_0} = Z \otimes V_1 \otimes V_1^{\otimes 2k}$
is also an $H(q)$-comodule. If $Z \otimes V_1^{\otimes N_0}$ was
a semisimple $H(q)$-comodule, it would be a semisimple
$\fre$-comodule, and 
$V_1^{\otimes N_0} = Z^{-1} \otimes Z \otimes V_1^{\otimes N_0}$
would be a semisimple $\fre$-comodule, and hence a semisimple 
$\qsl$-comodule. Thus $Z \otimes V_1^{\otimes N_0}$ is not a semisimple
$H(q)$-comodule: this concludes our proof. $\square$

\medskip
 
Proposition 4.2, combined with part a) of Theorem 1.1, proves 
part b) of Theorem 1.1.

\section{Corepresentations of $H(q)$, $q$ generic}

In this section $k$ is still an algebraically closed
field of characteristic zero, and $q \in k^*$ is generic.
We describe the simple $H(q)$-comodules and their
fusion rules, thereby completing the proof of Theorem 1.1.

\medskip

Let us begin with some preliminaries. We consider the monoid
$\N * \N$, the free product (=coproduct) of two copies of the monoid $\N$.
Equivalently $\N * \N$ is the free monoid on two generators
$\alpha$ and $\beta$ (this should not cause any confusion
with the elements $\alpha$ and $\beta$ of $H(q)$). 
There is a unique antimultiplicative morphism
$^- : \mathbb N * \mathbb N \longrightarrow \mathbb N * \mathbb N$
such that $\bar{e} = e$, $\bar{\alpha} = \beta$ and $\bar{\beta} = \alpha$
($e$ denotes the unit element of $\mathbb N * \mathbb N$).
Let $k[\N*\N]$ be the monoid algebra of $\N*\N$ : $k[\N*\N]$
is also the free algebra on two generators. Banica \cite{[Ba]} has
introduced a new product $\odot$ on $k[\N*\N]$. The following
lemma is Lemma 3 in \cite{[Ba]}, where the proof can be found.

\begin{lemm}
Consider the map $\odot : \N*\N\times \N*\N \longrightarrow
k[\N*\N]$ defined by
$$x\odot y = \sum_{x=ag,y = \bar{g}b} ab \ , \ x,y \in \N* \N \ ,$$
and extend $\odot$ to $k[\N*\N]$ by bilinearity. Then
$(k[\N*\N], + , \odot)$ is an associative $k$-algebra, with $e$ as unit
element. Furthermore $(k[\N*\N], + , \odot)$ is still the free algebra
on two generators: if $B$ is any algebra and $u,v \in B$, there exists
a unique algebra morphism \break
$\psi : (k[\N*\N], + , \odot) \longrightarrow B$
such that $\psi(\alpha) = u$ and $\psi(\beta) = v$. $\square$
\end{lemm}

We will need some character theory. Let $A$ be a Hopf algebra
and let
$V$ be a finite-dimensional $A$-comodule with corresponding coalgebra
map $\Phi_V : V^* \otimes V \longrightarrow A$.
Recall (see e.g. \cite{[KS]}) that the character of $V$ 
is defined to be $\chi_V := \Phi_V({\rm id}_V)$.
If $V$ and $W$ are  finite-dimensional $A$-comodules, then
$\chi(V \oplus W) = \chi(V) + \chi(W)$, 
$\chi(V \otimes W) = \chi(V) \chi(W)$
and $V \cong W \iff \chi(V) =  \chi(W)$.

\medskip

Recall \cite{[KP],[KS]} that $\qsl$ is cosemisimple and has a complete family
of simple comodules $(V_i)_{i \in \N}$, with $V_0=k$ and $\dim(V_i) = i+1$,
for $i \in \N$, and
$$V_i \otimes V_1 \cong V_1 \otimes V_i \cong V_{i-1} \oplus V_{i+1},
\ {\rm for} \ i \in \N^*. $$
As in the preceding section, for $i \in \Z$, we denote by
$Z^i$ the one-dimensional comodule corresponding to the element
$z^i$ of $k[z,z^{-1}]$. We identify $H(q)$ with a Hopf subalgebra
of $\fre$, via the morphism $\pi$ of Propositions 3.4 and 4.1.
Under this identification, the canonical two-dimensional comodules
$U$ and $V$ of $H(q)$ (see the notation in Section 1)
correspond to the simple alternated
comodules $Z \otimes V_1$ and $V_1 \otimes Z^{-1}$.

\begin{prop} There exists a unique algebra morphism
$\psi : (k[\N*\N],+ , \odot) \longrightarrow H(q)$ such that
$\psi(\alpha) = \chi(Z \otimes V_1)$ and $\psi(\beta)=
\chi(V_1 \otimes Z^{-1})$. Moreover for all $x \in \N * \N$, 
$\psi(x)$ is the character of a simple $H(q)$-comodule. 
\end{prop}

The first assertion is a direct consequence of Lemma 5.1. To prove the 
second one, we need a couple of lemmas. 

\begin{lemm} For all $n \in \N$, we have:
$$\psi((\alpha \beta)^n) = \chi(Z \otimes V_{2n} \otimes Z^{-1}) \quad ;
\quad   \psi((\beta \alpha)^n) = \chi(V_{2n}) \quad ;$$
$$\psi((\alpha \beta)^n \alpha ) = \chi(Z \otimes V_{2n+1})
\quad ; \quad \psi((\beta \alpha)^n\beta) = \chi(V_{2n+1} \otimes Z^{-1}).$$
\end{lemm}

\noindent
\textbf{Proof}. We prove the lemma by induction on $n$. For $n=0$,
the result is clear. Now assume that the lemma has been proved for 
$n\geq 0$. We have $(\alpha \beta)^n \alpha \odot \beta = 
(\alpha \beta)^{n+1} + (\alpha \beta)^n$, and so
\begin{align*}
\psi((\alpha \beta)^{n+1}) & = \psi((\alpha \beta)^n \alpha) \psi(\beta) 
- \psi((\alpha \beta)^n) \\
& = \chi(Z \otimes V_{2n+1}) \chi(V_1 \otimes Z^{-1}) - 
\chi(Z \otimes V_{2n} \otimes Z^{-1}) \ {(\rm by} \ {\rm induction)} \\
& = \chi(Z \otimes V_{2n} \otimes Z^{-1}) + 
\chi(Z \otimes V_{2n+2} \otimes Z^{-1}) - \chi(Z \otimes V_{2n} \otimes Z^{-1})
\\ 
& = \chi(Z \otimes V_{2(n+1)} \otimes Z^{-1}).
\end{align*}
Using $(\beta \alpha)^n \beta \odot \alpha = 
(\beta \alpha)^{n+1} + (\beta \alpha)^n$, one shows in the
same way that $\psi((\beta \alpha)^{n+1}) = \chi(V_{2(n+1)})$. 
We have $(\alpha \beta)^{n+1} \odot \alpha = 
(\alpha \beta)^{n+1}\alpha + (\alpha \beta)^n\alpha$, and hence
$$\psi((\alpha \beta)^{n+1}\alpha) = 
\psi((\alpha \beta)^{n+1}) \psi(\alpha) 
- \psi((\alpha \beta)^n\alpha).$$
We have already shown that $\psi((\alpha \beta)^{n+1}) =
\chi(Z \otimes V_{2(n+1)} \otimes Z^{-1})$, and by induction
$\psi((\alpha \beta)^n \alpha ) = \chi(Z \otimes V_{2n+1})$, so we have:
\begin{align*}
\psi((\alpha \beta)^{n+1}\alpha) & = \chi(Z  \otimes V_{2n+2} \otimes Z^{-1}
\otimes Z \otimes V_1) - \chi(Z \otimes V_{2n+1}) \\
& = \chi(Z \otimes V_{2n+1}) +\chi(Z \otimes V_{2n+3}) -
\chi(Z \otimes V_{2n+1}) = \chi(Z \otimes V_{2(n+1)+1}).
\end{align*}
One shows in a similar manner that
$\psi((\beta \alpha)^{n+1}\beta) = \chi(V_{2(n+1)+1} \otimes Z^{-1})$: this
concludes the proof. $\square$

\begin{lemm} Let $x \in \N * \N$. Then:

\noindent
$\bullet$ $\psi(x\alpha) = \chi(X \otimes V_i)$, for some $i \in \N^*$,
where $X=k$  or $X$ is a simple alternated comodule ending by $Z$ or
$Z^{-1}$.

\noindent
$\bullet$ $\psi(\alpha x) = \chi(Z \otimes X)$,
where $X$ is a simple alternated comodule beginning by some $V_i$, 
$i \in \N^*$. 

\noindent
$\bullet$ $\psi(x\beta) = \chi(X \otimes Z^{-1})$, where 
$X$ is a  simple alternated comodule ending by some $V_i$, 
$i \in \N^*$. 

\noindent
$\bullet$ $\psi(\beta x) = \chi(V_i \otimes X)$, for some $i \in \N^*$, 
where  $X = k$ or $X$ 
is a simple alternated comodule beginning by $Z$ or $Z^{-1}$. 
\end{lemm}

\noindent
\textbf{Proof}. We first prove the lemma for elements 
$x$ as in Lemma 5.3. Let $x = (\alpha \beta)^{n}$.
Then using Lemma 5.3, we have
\begin{align*}
& \psi  (x\alpha) = \psi((\alpha \beta)^{n}\alpha) = 
\chi(Z \otimes V_{2n+1}) \  , \\
& \psi(\alpha x) = \psi(\alpha(\alpha \beta)^n) = \psi(\alpha \odot
(\alpha \beta)^n)) = \psi(\alpha) \psi((\alpha \beta)^n) = \\
& \chi(Z \otimes V_1) \chi(Z \otimes V_{2n} \otimes Z^{-1}) = 
 \chi(Z \otimes V_1 \otimes Z \otimes V_{2n} \otimes Z^{-1}) \ , \\
& \psi  (x\beta) = \psi((\alpha \beta)^{n}\beta) =
 \chi(Z \otimes V_{2n} \otimes Z^{-1}) \chi(V_1 \otimes Z^{-1}) = 
\chi(Z \otimes V_{2n} \otimes Z^{-1} \otimes V_1 \otimes Z^{-1}) \ , \\
& \psi  (\beta x ) = \psi(\beta (\alpha \beta)^{n}) =
\psi((\beta \alpha)^n \beta) = \chi(V_{2n+1} \otimes Z^{-1}) \ .
\end{align*}
Similar computations show that the lemma is true for
$x= (\beta \alpha)^n$, $x= (\alpha \beta)^n\alpha$ or
$x = (\beta \alpha)^n \beta$.

We now prove the lemma for an arbitrary element $x \in \N * \N$ using 
an induction on the length $n$ of $x$. If $n= 0$, the result is obviously true.
Let us assume that the lemma has been proved for elements of length
$\leq n$ ($n \geq 0$), and let $x$ be an element of length $n+1$.
If $x$ is one of the elements of Lemma 5.3, the result has already been 
proved so we can assume that $x= y\alpha^2 z$ or that $x= y \beta^2 z$. 
For example assume that $x = y \alpha^2 z$. We have 
$$\psi(x\alpha) = \psi(y \alpha^2 z \alpha) = 
\psi(y \alpha \odot \alpha z \alpha) =
\psi(y \alpha) \psi(\alpha z \alpha).$$
By induction, we have $\psi(y \alpha) = X \otimes V_i$
for $i \in \N^*$ and $X=k$ or $X$ is a simple alternated
$\fre$-comodule ending by $Z$ or $Z^{-1}$. Also by induction
$ \psi(\alpha z \alpha) = \chi(Z \otimes Y \otimes V_j)$ for $j \in \N^*$, and 
$Y=k$ or $Y$ is a simple alternated comodule ending by $Z$ or $Z^{-1}$
and beginning by some $V_k$, $k \in \N^*$. So finally
$\psi(x\alpha) =  \chi(X \otimes V_i \otimes Z \otimes Y \otimes V_j)$
and $X \otimes V_i \otimes Z \otimes Y$ is a simple alternated comodule 
ending by $Z$ or $Z^{-1}$. We also have
$$\psi(\alpha x) = \psi(\alpha y \alpha^2 z) = 
\psi(\alpha y \alpha \odot z \alpha) =
\psi(\alpha y \alpha) \psi(\alpha z).$$
By induction we have $\psi(\alpha y \alpha)
= \chi(Z \otimes X \otimes V_i)$, $i \in \N^*$, and $X=k$ or 
$X$ is a simple alternated comodule beginning by some 
$V_j$, $j \in \N^*$ and ending by $Z$ or $Z^{-1}$. Also
$\psi(\alpha z) = \chi(Z \otimes Y)$ where $Y$ is a simple
alternated comodule beginning by some $V_k$, $k \in \N^*$.
Hence $\psi(\alpha x) = \chi(Z \otimes X \otimes V_i \otimes Z \otimes Y)$,
where $X \otimes V_i \otimes Z \otimes Y$ is a simple alternated comodule
beginning by some $V_j$, $j \in \N^*$. Let us now compute
$\psi(x \beta)$:
$$\psi(x\beta) = \psi(y \alpha^2 z \beta) = 
\psi(y \alpha \odot \alpha z \beta) =
\psi(y \alpha) \psi(\alpha z \beta).$$ 
By induction $\psi(y\alpha) = \chi(X \otimes V_i)$ where $X=k$
or $X$ is a simple alternated comodule ending by $Z$ or
$Z^{-1}$. Also $\psi(\alpha z \beta) = \chi(Z \otimes Y \otimes Z^{-1})$ where
$Y$ is a simple alternated comodule beginning by some $V_j$ and ending
by some $V_k$, $j,k \in \N^*$. So 
$\psi(x \beta) = \chi(X \otimes V_i \otimes Z \otimes Y \otimes Z^{-1})$,
where $X \otimes V_i \otimes Z \otimes Y$ is a simple 
alternated comodule ending by some $V_k$, $k \in \N^*$.
Let us finally compute $\psi(\beta x)$:
$$\psi(\beta x) = \psi(\beta y \alpha^2 z) = 
\psi(\beta y \alpha \odot \alpha z) =
\psi(\beta y \alpha) \psi(\alpha z).$$
By induction $\psi(\beta y \alpha) = \chi(V_i \otimes X \otimes V_j)$ for
$i,j \in \N^*$, and $X$ is a simple alternated comodule 
beginning by $Z$ or $Z^{-1}$ and ending by $Z$ or $Z^{-1}$.
Also $\psi(\alpha z) = \chi(Z \otimes Y)$, where $Y$ is an
alternated simple comodule beginning by some $V_k$, $k \in \N^*$.
So $\psi(\beta x) = \chi(V_i \otimes X \otimes V_j \otimes Z \otimes Y)$
where $X \otimes V_j \otimes Z \otimes Y$ is a simple alternated comodule
beginning by $Z$ or $Z^{-1}$.
Very similar computations prove the result for $x= y \beta^2 z$, and conclude
the proof of Lemma 5.4. $\square$ 

\medskip

Proposition 5.2 is a direct consequence of Lemma 5.4. $\square$

\bigskip

We can now easily list the simple $H(q)$-comodules,
and describe their fusion rules.
For $x \in \N * \N$, let $U_x$ be a simple $H(q)$-comodule
such that $\chi(U_x) = \psi(x)$. We have $U_e=k$, $U_\alpha = U$
and $U_\beta = V$, for the notations of the introduction. We have
$$\chi(U_x \otimes U_y) = \chi(U_x) \chi(U_y) = \psi(x) \psi(y)
= \psi(x\odot y) = \psi(\sum_{x=ag,y=\bar{g}b}ab)
= \chi(\bigoplus_{x=ag,y = \bar{g}b}U_{ab}),$$
and hence
$$U_x \otimes U_y \cong \bigoplus_{x=ag,y = \bar{g}b} U_{ab} \ .$$
By Lemma 5.4 we have $U_x \cong k$ if and only if $x = e$, and 
using the last formula, we see that Hom$(k, U_x \otimes U_y) \not = (0)$
if and only if $y= \bar{x}$. This implies that $U_x^* \cong U_{\bar{x}}$
and that $U_x  \cong U_y$ if and only if $x=y$.
Thus we have a family of simple $H(q)$-comodules
$(U_x)_{x \in \N * \N}$ whose coefficients generate $A$ as an
algebra, containing the trivial comodule and stable under tensor products:
using e.g. the orthogonality relations \cite{[KS]} we conclude that
any simple $H(q)$-comodule is isomorphic with a comodule $U_x$.

\medskip

The preceding discussion concludes the proof of Theorem 1.1:
there just remain to be said that the
monoidal category equivalence Comod$(H(F)) \cong^\otimes {\rm Comod}(H(q))$
transforms the fundamental $n$-dimensional comodules $U$ and $V$ of
$H(F)$ into the fundamental $2$-dimensional comodules
$U$ and $V$ of $H(q)$.

\medskip

Lemma 5.1, Proposition 5.2 and Theorem 1.1 combined together also 
yield the description of the Grothendieck $K_0$-ring of the
category Comod$_{\rm f}(H(F))$.

\begin{coro} Let $F \in GL(n,k)$ ($n \geq 2$) be a generic matrix. Then
we have a ring isomorphism
$$K_0({\rm Comod}_{\rm f}(H(F)) \cong \Z\{X,Y\} \ .$$ 
\end{coro}  

\section{Some applications}

We use Theorem 1.1 to prove a few structural results concerning
the Hopf algebras $H(F)$, for generic matrices.
Again $k$ is an algebraically closed field of characteristic zero.

\smallskip

Let us begin with the isomorphic classification. 
For universal compact quantum groups, this was done by Wang \cite{[Wa3]}.
Since we use the same type of arguments, we
will be a little concise.

\begin{prop}
Let $E \in GL(m,k)$, $F \in GL(n,k)$ ($m,n \geq 2$)
be generic matrices. The Hopf algebras
$H(E)$ and $H(F)$ are isomorphic if and only if one of the
two conditions hold.

\noindent
i) $m=n$ and there exists $P \in GL(n,k)$ such that $F = \pm PEP^{-1}$.

\noindent
ii) $m=n$ and there exists $P \in GL(n,k)$ such that 
${^t \! F}^{-1} = \pm P E P^{-1}$.
\end{prop}

\noindent
\textbf{Proof}. Let $f : H(E)  \longrightarrow H(F)$ be a Hopf algebra 
isomorphism, and denote by $f_* : {\rm Comod}(H(E)) \longrightarrow
{\rm Comod}(H(F))$ the functor induced by $f$. By \cite{[Wa3]} $U$
and $V$ are the simple $H(E)$-comodules (resp. $H(F)$-comodules) 
with the strictly smallest dimension, and hence we have
$f_*(U) \cong U$ or $f_*(U) \cong V$. If $f_*(U) \cong U$, then
$m=n$ and there exists $P \in GL(n,k)$ such that 
$f(u) = {^t \! P} u {^t \! P}^{-1}$ and necessarily
$f(v) = P^{-1}vP$. Since $f$ is well-defined and since $U$ is simple, 
it is easy to check that $F = \pm PEP^{-1}$.  
If $f_*(U) \cong V$, then
$m=n$ and there exists $P \in GL(n,k)$ such that 
$f(u) = {^t \! P} v {^t \! P}^{-1}$ and necessarily
$f(v) = P^{-1}{^t \! F}^{-1} v {^t \! F}P$. 
Since $f$ is well-defined and since $U$  and $V$ are simple, 
it is easy to check that ${^t \! F}^{-1} = \pm PEP^{-1}$. 

Conversely, if  $F = \pm PEP^{-1}$, it is easy to check that
there exists a Hopf algebra isomorphism
$f : H(E) \longrightarrow H(F)$ such that 
$f(u) = {^t \! P} u {^t \! P}^{-1}$ and 
$f(v) = P^{-1}vP$. 
if  ${^t \!F}^{-1} = \pm PEP^{-1}$, it is easy to check that
there exists a Hopf algebra isomorphism
$f : H(E) \longrightarrow H(F)$ such that 
$f(u) = {^t \! P} v {^t \! P}^{-1}$ and 
$f(v) = P^{-1}{^t \! F}^{-1}u {^t \! F}P$. $\square$

\bigskip

Let us now compute the automorphism group of the Hopf algebra $H(F)$.
Let $F \in GL(n,k)$. Put
$$X_0(F) = \{K \in GL(n,k) \ | \ KFK^{-1} = F \} \ , \
Y(F) = \{K \in GL(n,k) \ | \ KFK^{-1} = {^t \! F}^{-1} \},$$
and $X(F) = X_0(F)/k^*$. 
Then $X(F)$ is a group. For $N \in \N \cup \{\infty\}$,
the cyclic group of order $N$ is denoted by $C_N$.

\begin{prop}
Let $F \in GL(n,k)$ ($n \geq 2$) be a generic matrix.

\noindent
a) Assume that $Y(F) = \emptyset$. Then $X(F) \cong
{\rm Aut}_{\rm Hopf}(H(F))$.

\noindent
b) Assume that $Y(F) \not = \emptyset$. 
Let $K \in Y(F)$, and put
$N = {\rm min}\{p \in  \N \cup \{\infty\} \ | \ 
(F^{-1} {^t \!K}^{-1}K)^p \in k^*
\}$. Then we have an exact sequence of groups
$$
1 \longrightarrow C_N \longrightarrow 
X(F) \rtimes C_{2N} \longrightarrow {\rm Aut}_{\rm Hopf}(H(F))
\longrightarrow 1$$
In particular, if there exists $K \in GL(n,k)$ such that
$F= {^t \! K}^{-1}K$, then $K \in Y(F)$, we can take $N=1$ and we
have an isomorphism $X(F) \rtimes C_{2} \cong {\rm Aut}_{\rm Hopf}(H(F))$.
\end{prop}

\noindent
\textbf{Proof}. Let $K \in X_0(F)$. Then
there exists a Hopf algebra automorphism $\phi_K$ of
$H(F)$ such that  
$\phi_K(u) = {^t \! K} u {^t \! K}^{-1}$ and 
$\phi_K(v) = K^{-1}vK$. This gives a group morphism
$\phi : X(F) \longrightarrow  {\rm Aut}_{\rm Hopf}(H(F))$, injective since
the comodule $U$ is simple. Now consider $f \in {\rm Aut}_{\rm Hopf}(H(F))$.
Then by the proof of Proposition 6.1, either there exists
$K \in X_0(F)$ (recall that tr$(F) \not = 0$ since $F$ is generic)
such that
$f(u) = {^t \! K} u {^t \! K}^{-1}$, either there exists
$K \in Y(F)$ such that $f(u) =  {^t \! K} v {^t \! K}^{-1}$.
If $Y(F) = \emptyset$, then $f= \phi_K$ and the morphism 
$\phi$ is an isomorphism.
Assume now that $Y(F) \not = \emptyset$ and let $K \in Y(F)$.
Then there exists $\psi_K \in   {\rm Aut}_{\rm Hopf}(H(F))$ such
that $\psi_K(u) =  {^t \! K} v {^t \! K}^{-1}$ and 
$\psi_K(v) = K^{-1}{^t \! F}^{-1}u {^t \! F}K$. Let
$f \in {\rm Aut}_{\rm Hopf}(H(F))$. Then by the proof of Proposition
6.1 there exists $K \in X_0(F)$ such that $f = \phi_K$ or there
exists $M \in Y(F)$ such that $f = \psi_M$. 
We have $\psi_M = \phi_{{^t \! M}^{-1} {^t \! K}} \circ \psi_K$, and 
thus $G = \phi(X(F))\langle \psi_K \rangle$. For $M,L \in Y(F)$, we have
$\psi_M  \circ \psi_L = \phi_{F^{-1} {^t \! M}^{-1} L}$, hence
$|\langle \psi_K \rangle| = 2N$. Also
$\langle \psi_K \rangle \cap \phi(X(F)) = \langle 
\phi_{F^{-1}{^t \! K}^{-1}K} \rangle$ and   
$|\langle \psi_K \rangle \cap \phi(X(F))| = N$. We have 
$\psi_K \circ \phi_L \circ \psi_K^{-1} 
= \phi_{F^{-1}{^t \! K}^{-1} {^t \! L}^{-1} {^t \! K}F}$ and hence
$\phi(X(F))$ is a normal subgroup of ${\rm Aut}_{\rm Hopf}(H(F))$. 
We can now use a well-know result in group theory:
if $G$ is a group with two subgroups $H$ and $ K$ such that
$G=HK$, such that $H$ is normal in $G$ and such that
$H \cap K$ is abelian, then we have a group exact sequence
$$
1 \longrightarrow H \cap K \longrightarrow H \rtimes K
\longrightarrow G \longrightarrow 1.$$
The last assertion is immediate. $\square$

\section{Quantum automorphism groups of matrix algebras} 

In his paper \cite{[Wa2]}, Wang described the quantum automorphism 
group of a finite-dimensional $C^*$-algebra endowed with a trace,
the term quantum automorphism group (or quantum symmetry group)
being understood in the sense of Manin \cite{[Ma]}. We refer
the reader to \cite{[Ma]} or \cite{[Wa2]} for these ideas.
The representation theory of such quantum automorphism
groups was described by Banica \cite{[Ba2]} in the case of good traces, and
is similar to the one of $SO(3)$.

In \cite{[Bi1]} we proposed a natural categorical generalization 
of Wang's construction, yielding in particular an algebraic
analogue of the quantum automorphism group of a finite-dimensional
measured algebra.
We will see that in the case of a measured matrix algebra
with a non-necessarily tracial measure,
the results of the present paper enable us to describe 
the representation theory of such a quantum group, reducing
the computations to the case of the quantum $SO(3)$-group.

\medskip

Recall \cite{[Bi1]} that a measured algebra is a pair $(Z, \phi )$
where $Z$ is an algebra and $\phi : Z \longrightarrow k$ is a linear
map such that the bilinear form $Z \times Z \longrightarrow k$,
$(a,b) \mapsto \phi(ab)$, is
non-degenerate. 
We will only be concerned here by the example $(M_n(k), {\rm tr_F})$
where $F \in GL(n,k)$ and ${\rm tr}_F = {\rm tr}({^t \! F}^{-1}-)$. The quantum
automorphism groups of  $(M_n(k), {\rm tr_F})$, denoted
$A_{aut}(M_n(k), {\rm tr_F})$, may be described as follows
(see \cite{[Wa2]} for details).
As an algebra $A_{aut}(M_n(k), {\rm tr_F})$ is the universal algebra
with generators $X_{ij}^{kl}$, $1 \leq i,j,k,l \leq n$, and satisfying 
the relations ($1 \leq i,j,k,l,r,s \leq n$):
$$
\sum_t X_{rt}^{ij} X_{ts}^{kl} = \delta_{jk} X_{rs}^{il} \ ;
\ \sum_{t,p}F_{tp} X_{kl}^{it} X_{rs}^{pj} = F_{lr}X_{ks}^{ij} \ ; \
\sum_t X_{ij}^{tt} = \delta_{ij} \ ; \ 
\sum_{t,p}F_{tp}^{-1}X_{tp}^{ij} = F_{ij}^{-1}.$$
It has a natural Hopf algebra structure given by
$$\Delta(X_{ij}^{kl}) = \sum_{r,s} X_{ij}^{rs} \otimes X_{rs}^{kl} \ ; \
\varepsilon(X_{ij}^{kl}) = \delta_{ik} \delta_{jl} \ ; \
S(X_{ij}^{kl}) = \sum_{r,s} F_{jr}F_{sl}^{-1} X_{sk}^{ri} \ , 
1 \leq i,j,k,l \leq n.$$
Let $E \in GL(m,k)$ and $F \in GL(n,k)$. Let us define
the algebra  {\small$A_{is}(M_m(k), {\rm tr_E}; M_n(k), {\rm tr}_F)$}
to be the universal algebra
with generators $X_{ij}^{kl}$, $1 \leq i,j\leq m$, $1 \leq k,l \leq n$, 
and satisfying 
the relations:
\begin{align*}
& \sum_{t=1}^m  X_{rt}^{ij} X_{ts}^{kl} = \delta_{jk} X_{rs}^{il} \quad ,
 \quad 1 \leq i,j ,k,l \leq n \quad  , \quad 1 \leq r,s \leq m \ ;\\
&  \sum_{t,p=1}^nF_{tp} X_{kl}^{it} X_{rs}^{pj} = E_{lr}X_{ks}^{ij} \quad ,
\quad 1 \leq i,j \leq n \quad , \quad 1 \leq k,l,r,s \leq m \ ; \\
 & \sum_{t=1}^n X_{ij}^{tt} = \delta_{ij}\ , \ 1 \leq i,j \leq m \quad ; \quad 
\sum_{t,p=1}^m E_{tp}^{-1}X_{tp}^{ij} = F_{ij}^{-1} \quad , \quad
1 \leq i,j \leq n.
\end{align*}

\begin{lemm}
Let $E \in GL(m,k)$ and let $F \in GL(n,k)$ ($m,n \geq 2$) with 
${\rm tr}(E) = {\rm tr}(F)$ and ${\rm tr}(E^{-1}) = {\rm tr}(F^{-1})$. 
Then $A_{is}(M_m(k), {\rm tr}_E ; M_n(k), {\rm tr}_F)$ is a non-zero algebra.
\end{lemm}

\noindent
\textbf{Proof}. It is straightforward to check that there exists a unique 
algebra morphism $\varphi : A_{is}(M_m(k), {\rm tr}_E ; M_n(k), {\rm tr}_F)
\longrightarrow H(E,F)$ such that $\varphi(X_{ij}^{kl}) = u_{ik}v_{jl}$
for $1 \leq i,j \leq n$, $1 \leq k,l \leq m$. The elements 
$u_{ik}v_{jl}$ are non-zero elements of $H(E,F)$ by Section 2, and hence
$A_{is}(M_m(k), {\rm tr}_E ; M_n(k), {\rm tr}_F)$ is a non-zero algebra.
$\square$ 

\bigskip

We arrive at the main result of the section:

\begin{theo}
Let $E \in GL(m,k)$ and let $F \in GL(n,k)$ ($m,n \geq 2$) with 
${\rm tr}(E) = {\rm tr}(F)$ and ${\rm tr}(E^{-1}) = {\rm tr}(F^{-1})$. 
Then the comodule categories over $A_{aut}(M_m(k),{\rm tr}_E)$ and
$A_{aut}(M_n(k),{\rm tr}_F)$ are monoidally equivalent.
In particular, if ${\rm tr}(F) = {\rm tr}(F^{-1})$
and if there exists $q \in k^*$ such that $q^2 -{\rm tr}(F)q +1 = 0$, then
the comodule categories over $A_{aut}(M_n(k),{\rm tr}_F)$ and 
$\mathcal O (SO_{q^{1/2}}(3))$ are monoidally equivalent.
\end{theo} 

\noindent
\textbf{Proof}. Let us show that {\small 
$$(A_{aut}(M_m(k),{\rm tr}_E), A_{aut}(M_n(k),{\rm tr}_F),
A_{is}(M_m(k),{\rm tr}_E; M_n(k),{\rm tr}_F),
A_{is}(M_n(k),{\rm tr}_F; M_m(k), {\rm tr}_E))$$} is a Hopf-Galois
system \cite{[Bi4]}. First by Lemma 7.1 all these algebras are non-zero.
Let $G \in GL(p,k)$. It is a direct computation to check that there
exists a unique algebra morphism {\small
$$\delta_{E,F}^G : A_{is}(M_m(k),{\rm tr}_E; M_n(k),{\rm tr}_F)
\longrightarrow A_{is}(M_m(k),{\rm tr}_E; M_p(k),{\rm tr}_G)
\otimes A_{is}(M_p(k),{\rm tr}_G; M_n(k),{\rm tr}_F)$$}
such that $\delta_{E,F}^G(X_{ij}^{kl}) = \sum_{r,s}
X_{ij}^{rs} \otimes X_{rs}^{kl}$. Also there exists a unique
algebra morphism $\phi : A_{is}(M_n(k),{\rm tr}_F;M_m(k),{\rm tr}_E)
\longrightarrow 
A_{is}(M_m(k),{\rm tr}_E; M_n(k),{\rm tr}_F)^{\rm op}$ such that
$\phi(X_{ij}^{kl}) = \sum_{r,s} F_{jl}E_{sl}^{-1}X_{sk}^{ri}$.
With these structural morphisms, it is immediate to check that we
indeed have a Hopf-Galois system. Hence using Corollary 1.4 of 
\cite{[Bi4]}, we have our monoidal category equivalence.
Now assume that 
${\rm tr}(F) = {\rm tr}(F^{-1})$
and that there exists $q \in k^*$ such that $q^2 -{\rm tr}(F)q +1 = 0$.
Put ${\rm tr}_q = {\rm tr}_{F_q}$. Then we have an
equivalence of monoidal categories:
$$A_{aut}(M_n(k),{\rm tr}_F) \cong^\otimes 
A_{aut}(M_2(k), {\rm tr}_q).$$ 
Finally it may be shown that $A_{aut}(M_2(k), {\rm tr}_q)$ and
$\mathcal O (SO_{q^{1/2}}(3))$ are isomorphic. One considers first
the Hopf algebra morphism $A_{aut}(M_2(k), {\rm tr}_q) \longrightarrow
\qsl$ obtained using the adjoint corepresentation of the canonical
two-dimensional $\qsl$-comodule.
This Hopf algebra morphism is injective, and using
\cite{[Di]}, we arrive at the desired conclusion.
$\square$

\end{document}